\documentclass[11pt, a4paper, oneside, final, reqno]{amsart}

\usepackage{amssymb}
\usepackage{comment}

\newcommand{\calA}{\mathcal{A}}
\newcommand{\calB}{\mathcal{B}}

\newcommand{\calX}{\mathcal{X}}
\newcommand{\calY}{\mathcal{Y}}

\newcommand{\bbN}{\mathbb{N}}
\newcommand{\bbZ}{\mathbb{Z}}

\newcommand{\bbR}{\mathbb{R}}
\newcommand{\bbC}{\mathbb{C}}

\newcommand{\bbF}{\mathbb{F}}
\newcommand{\bbL}{\mathbb{L}}

\newcommand{\id}{\operatorname{id}}

\newcommand{\mk}{\operatorname{mk}}

\newcommand{\tens}{{\otimes}}

\usepackage{thmtools}
\declaretheorem[style=theorem,name={Theorem}]{theoremletter}

\theoremstyle{definition}
\newtheorem{definition}{Definition}[section]
\newtheorem{notation}{Notation}[section]

\theoremstyle{plain}
\newtheorem{theorem}[definition]{Theorem}
\newtheorem{proposition}[definition]{Proposition}
\newtheorem{corollary}[definition]{Corollary}
\newtheorem{lemma}[definition]{Lemma}

\theoremstyle{remark}
\newtheorem{remark}[definition]{Remark}

\numberwithin{equation}{section}
\usepackage{hyperref}
\hypersetup{pdftitle={}}
\usepackage[nameinlink, noabbrev, capitalize]{cleveref}
\usepackage[utf8]{inputenc}

\usepackage{xcolor}
\usepackage[all]{xy}

\usepackage[a4paper, total={6.5in, 10in}]{geometry}
\usepackage{microtype}

\usepackage{todonotes}

\title[]{Cartesian products of groups as quantum metric spaces}

\date{\today}

\author[1]{Are Austad}
\address{Are Austad, Department of Mathematics, University of Oslo, P.O.Box 1053 Blindern, 0316 Oslo, Norway}
\email{areaus@math.uio.no}

\author[2]{David Kyed}
\address{David Kyed, Department of Mathematics and Computer Science, University of Southern Denmark, Campusvej 55, DK-5230 Odense M, Denmark}
\email{dkyed@imada.sdu.dk}

\begin{document}
	
	\maketitle
	
	\begin{abstract}
        We provide a new class of finitely generated groups whose word-length functions give rise to compact quantum metric spaces. Kaad's recent theory of matrix compact quantum metric spaces is instrumental in obtaining this new class,  and along the way we also show that any length function on an amenable group which gives rise to a compact quantum metric space in fact gives rise to a matrix compact quantum metric space, and that the corresponding result is false for the free group on two generators.

	\end{abstract}

	\section{Introduction}

Over the past decade, substantial progress has been made in developing Rieffel's theory of compact quantum metric spaces, which provides a far-reaching non-commutative generalisation of classical compact metric spaces, with strong ties to Connes' non-commutative geometry and theoretical physics. Recent advances in the field \cite{AKK:Podcon, Leimbach2025, walter:GH-convergence, walter-connes:truncations}, suggest that basing the theory on operator systems, rather than order unit spaces or full $C^*$-algebras, provides a good blend of flexibility and rigidity, and with this in mind it is also very natural to seek a matricial version of the theory. A first  approach to a matricial theory was developed by Kaad in \cite{Kaad24}, and the resulting theory of matrix compact quantum metric spaces has, among other things, the advantage that it admits a natural tensor product construction, compatible with the external product of spectral triples. The present paper continues the investigation of matrix compact quantum metric spaces, with particular emphasis on the examples arising from groups. These examples fit into the more general connection between compact quantum metric spaces and non-commutative (differential) geometry, which we now describe in more detail. Firstly, a compact quantum metric space consists of an operator system $\calX$ together with a seminorm $L\colon \calX \to [0,\infty)$, such that the associated  {Monge-Kantorovi\v{c} metric} 
\begin{align*}
		\operatorname{mk_L}(\phi, \psi) := \sup \big\{ \vert \phi(x) - \psi(x) \vert \mid x \in \calX, L(x) \leq 1 \big\}, \qquad \phi,\psi \in S(\calX),
	\end{align*}
	induces the weak$^*$ topology on the state space $S(\calX)$. 
The seminorm $L$  may be viewed as a generalisation of the example arising from a  classical compact metric space $(M,d) $, where $\calX$ is  the algebra of Lipschitz continuous functions and $L$ the seminorm that assigns the (smallest) Lipschitz constant to the function in question. In the same vein, a matrix compact quantum metric space consists of an operator system $\calX$ together with a family of seminorms $L_n\colon M_n(\calX) \to [0,\infty)$ satisfying a natural set of axioms; see Section \ref{sec:prelims} for details. \\
On the other hand, the fundamental object in Connes' non-commutative geometry \cite{Connes1994} is the notion of a spectral triple, and from a spectral triple $(\calA, H,D)$ one obtains a natural seminorm $L_D\colon \calA \to [0,\infty)$ given by $L_D(a):=\|[D,a]\|$. In many interesting examples, but not always, this gives rise to a compact quantum metric space.  From a discrete group $G$ equipped with a proper length function $\ell\colon G\to [0,\infty)$,  one obtains a natural spectral triple $(\bbC[G], \ell^2(G), D_\ell)$, and even in this situation it has proven very difficult to determine when one obtains a compact quantum metric structure from the associated seminorm $L_{D_\ell}$. On the positive side, this has been verified for hyperbolic groups and groups of polynomial growth in \cite{OzawaRieffel} and \cite{ChristRieffel}, respectively, and only very recently the first counterexample for a word-length function was found by Klisse \cite{Klisse2026}.  \\
The aim of the present paper is two-fold: First, we answer,  negatively,  a question posed in \cite{Kaad24}, asking if any compact quantum metric space is in fact automatically a matrix compact quantum metric space. Kaad conjectured that this is not the case, and we confirm this conjecture by showing that the free group on two generators provides a natural counterexample. 

\begin{theoremletter}[see Theorem \ref{thm:free-group-not-matrix}]\label{thm:mainA}
For the free group $\bbF_2$ on two generators, $a$ and $b$, the word-length function arising from the generating set $S:=\{a,a^{-1}, b,b^{-1}\}$ gives rise to a compact quantum metric space, which is not a matrix compact quantum metric space.
\end{theoremletter}

Contrasting this, we provide an affirmative answer in the presence of amenability.
\begin{theoremletter}[see Theorem \ref{thm:amenable-CQMS-is-MCQMS}]\label{thm:mainB}
If $G$ is an amenable group equipped with a proper length function yielding a compact quantum metric space, then it is automatically a matrix compact quantum metric space. 
\end{theoremletter}

Secondly, exploiting the tensor product construction for matrix compact quantum metric spaces, we provide a new class of groups for which one does obtain a compact quantum metric structure.

\begin{theoremletter}[see Theorem \ref{thm:cartesian-products}]\label{thm:mainC}
Let $G$ and $K$ be finitely generated groups equipped with word-length functions $\ell_G$ and $\ell_K$. If $(\bbC[G], L_{D_{\ell_G}})$ is a compact quantum metric space and $K$ is either finite, abelian or the infinite dihedral group, then $(\bbC[G\times K], L_{D_\ell})$ is a compact quantum metric space for the length function $\ell(g,k)=\ell_G(g)+\ell_K(k)$.
\end{theoremletter}

Theorem \ref{thm:mainC} does indeed provide genuinely new examples, in that one may, for instance, take $G$ to be any non-elementary hyperbolic group and $K$ to be the integers, then $G\times K$ is neither hyperbolic (as it contains a copy of $\bbZ^2$) nor of polynomial growth (as $G$ must have exponential growth), and is thus not covered by the existing results in \cite{OzawaRieffel} and \cite{ChristRieffel}.

	\subsubsection*{Notation} 
Throughout the article, all $C^*$-algebras are unital, and all spectral triples are faithfully and unitally represented. Moreover, inner products are assumed linear in the first variable. 
	Algebraic tensor products are denoted by $\odot$, while $\otimes_{\rm min}$ will denote the minimal tensor product of $C^*$-algebras, while $\otimes$ will be used for tensor products in the category of Hilbert spaces. A unital completely positive map between operator systems will be referred to as a UCP map.

	\subsubsection*{Acknowledgments} 
	The present paper grew out of a series of interesting discussions with Jens Kaad and Nadia S.~Larsen in the fall of 2025, and it is a pleasure to thank them both for sharing their insights. We also wish to acknowledge that the fundamental idea of geodesic cancellation, which plays a pivotal role in the proof of Theorem \ref{thm:mainC}, was originally discovered by Jens  in the case of the integers.  That Theorem \ref{thm:amenable-CQMS-is-MCQMS}  and Theorem \ref{thm:one-MCQMS-suffices} should be true, was also suggested by Jens.\\
The first named author would also like to thank the University of Southern Denmark for its hospitality during a visit in November 2025, during which some of this work was conducted. The first named author was funded  by the Research Council of Norway projects 324944 and 359602. The second named author  was funded by  the Independent Research Fund Denmark through grant no.~1026-00371B and by  the European Commission through the MSCA Staff Exchanges grant no.~101086394. 

\subsubsection*{AI declaration} 
The present paper consists of a mixture of human-generated results and AI-assisted proofs. As indicated above, the fundamental idea of geodesic cancellation leading to Theorem \ref{thm:mainC}  was discovered and utilised by the authors (without AI assistance but through conversations with Jens Kaad and Nadia S.~Larsen) in the case where $K$ is $\bbZ$ or $\bbZ^2$, but the full class of groups covered in Theorem \ref{thm:mainC} was then obtained via ChatGPT through a carefully designed research prompt. Similarly, Theorem \ref{thm:mainB} was obtained without any AI assistance, whereas the counterexample in Theorem \ref{thm:mainA} was obtained through a dialogue with ChatGPT,  to explore if Theorem \ref{thm:mainB}  could be extended to also cover the class of hyperbolic groups (which is then shown to be impossible by Theorem \ref{thm:mainA}).

	\section{Preliminaries}
	\label{sec:prelims}
	
	The present section introduces the necessary background concerning compact quantum metric spaces and their matrix analogues. In the interest of brevity, we  will start by introducing the matrix version.
	
	\subsection{Matrix compact quantum metric spaces}
	Following \cite{Kaad24}, we choose the operator system approach to  (matrix) compact quantum metric spaces. By an operator system, we mean a unital, $*$-invariant subspace $\mathcal{X}$ of a unital $C^*$-algebra $A$. Note that we do not require it to be norm closed. As defined, the ambient $C^*$-algebra is part of the definition of the operator system $\mathcal{X}$, but we remark that the theory also admits an abstract description; see \cite{PaulsenCBOA02} for details.

	We will tacitly identify the scalars inside an operator system $\mathcal{X}$ as $\bbC \cdot 1_{\mathcal X}$.  
A \emph{Lipschitz seminorm} on an operator system $\mathcal{X}$,  is a $*$-invariant seminorm $L\colon \operatorname{Dom}(L)\to [0,\infty)$, defined on a dense operator-subsystem $\operatorname{Dom}(L)\subseteq \mathcal{X}$, which vanishes on the unit.  As is common in the field, we will sometimes regard such a seminorm as an extended seminorm on the whole operator system by setting it to $\infty$ outside its domain.
	To discuss matrix compact quantum metric spaces, we will adopt the terminology developed in \cite{Kaad24}. 
	A \emph{Lipschitz operator seminorm} on an operator system $\calX$ is a family $\bbL = (L_n)_{n\geq 1}$ of $*$-invariant seminorms on $M_n(\calX)$, which vanish on $M_n(\bbC)$, such that
	\begin{align*}
		L_{n+m}(x \oplus y) \leq \max \{ L_n(x) , L_m(y) \} \qquad \text{ and } \qquad L_n(\alpha x \beta) \leq \Vert \alpha \Vert \cdot L_n(x)\cdot \Vert  \beta \Vert, 
	\end{align*}
	for all $x\in M_n(\calX)$, $y\in M_m(\calX)$ and scalar matrices $\alpha, \beta \in M_n(\bbC)$. \\
	
	\begin{notation}
	
	For an operator system $\calX$ and an $n\in \mathbb{N}$, we will denote by  $[\ \cdot \ ] \colon M_n (\calX) \to M_n(\calX)/M_n(\bbC)$  the natural quotient map. Moreover,  for a linear map $T\colon \calX \to \calY$ between operator systems,  we denote by $T_n\colon M_n(\calX)\to M_n(\calY)$ its entrywise matrix amplification.
	\end{notation}
	\begin{definition}[{{\cite{Kaad24}}}]\label{def:MCQMS}
		A pair $(\calX, \bbL)$ consisting of an operator system $\calX$ and a Lipschitz operator seminorm $\bbL=(L_n)_{n\in \bbN}$ is said to be a \emph{matrix compact quantum metric space} (MCQMS) if the following conditions hold:
		\begin{enumerate}
			\item There exists an $R>0$ such that
			\begin{align*}
				\Vert [x] \Vert \leq R\cdot  L_n(x), \qquad \text{for all }n\geq 1, \, x \in M_n(\calX).
			\end{align*}
			\item For every $\varepsilon > 0$ there exist an operator system $\mathcal{Y}_\varepsilon$, a unital complete isometry $\iota_\varepsilon \colon \calX \to \mathcal{Y}_\varepsilon$, and a finite-rank unital completely positive map $\Phi_\varepsilon\colon \calX \to \mathcal{Y}_\varepsilon$ satisfying
			\begin{align*}
				\Vert (\iota_\varepsilon)_n(x) - (\Phi_\varepsilon)_n(x) \Vert \leq \varepsilon \cdot L_n(x) \qquad \text{for all } n\geq 1, \, x \in M_n(\calX). 
			\end{align*}
		\end{enumerate}
		Note in particular that $R$ and the approximating maps are independent of the matrix size $n$. For a single Lipschitz seminorm $L$, $(\calX, L)$ is called a \emph{compact quantum metric space} (CQMS) if these conditions hold at $n=1$, with $L_1 = L$, and where $\iota_\varepsilon$ is only required to be a unital linear isometry and $\Phi_\varepsilon$ a bounded, unital positive finite-rank linear map. In this case,  $L$ is called a \emph{Lip-norm}.
	\end{definition}

	Note first that any MCQMS is a CQMS by restricting the matrix approximations to the $n=1$ case, see also \cite[Proposition 5.4]{Kaad24}. Furthermore, the notion of CQMS from Definition \ref{def:MCQMS} agrees with Rieffel's usual definition by \cite[Theorem 3.1]{Kaad24}. In particular,  a Lipschitz seminorm $L \colon \calX \to [0,\infty)$ is a Lip-norm, if and only if the \emph{Monge-Kantorovi\v{c} metric} 
	\begin{align*}
		\operatorname{mk_L}(\phi, \psi) := \sup \{ \vert \phi(x) - \psi(x) \vert \mid x \in \calX, L(x) \leq 1 \}, \qquad \phi,\psi \in S(\calX),
	\end{align*}
	induces the weak$^*$ topology on the state space $S(\calX)$. Given a classical compact metric space $(M,d)$,
	one obtains \cite[Proposition 5.7]{Kaad24} an MCQMS by setting $\operatorname{Dom}(L_n):=\{f\colon M\to M_n(\bbC) \mid f \text{ Lipschitz continuous}\}$ and
\[
L_n(f):=\sup\left\{\frac{\|f(p)-f(q)\|}{d(p,q)} \mid p,q\in M, p\neq q\right\}
\]	 
Here we tacitly identify $C(M,M_n(\bbC))$ with $M_n(C(M))$.\\

We end this section, by recalling the following well known useful criterion.
	\begin{proposition}\label{prop:totally-bounded-criterion}
		Fix $\mu \in S(\calX)$, and put 
		\begin{align*}
			E_{L, \mu} = \big\{x \in \calX \mid L(x) \leq 1, \mu(x) = 0\big\}.
		\end{align*}
		Then $(\calX, L)$ is a CQMS if and only if $E_{L, \mu}$ is totally bounded in norm. In this case, $\ker L = \bbC$, the radius $R_\mu := \sup_{x \in E_{L,\mu}} \Vert x \Vert$ is finite, and
		\begin{align*}
			\Vert x - \mu(x) 1 \Vert \leq R_\mu\cdot  L(x)
		\end{align*}
		for all $x \in \calX$. Moreover,  for every $s,t \geq 0$, the set $\{x \in \calX \mid L(x) \leq s, \vert \mu(x) \vert \leq t \}$ is totally bounded.  
	\end{proposition}
	\begin{proof}
		The statement on total boundedness of $E_{L,\mu}$ follows from \cite[Proposition 1.3]{OzawaRieffel}. It is well known that if $(\calX, L)$ is a CQMS, then $\ker L = \bbC$; see e.g.~\cite[Lemma 2.2]{KaadKyed2021} for a proof. The statement $\Vert x - \mu(x) 1 \Vert \leq R_\mu L(x)$ follows from noting that $x - \mu(x) 1 \in L(x) \cdot E_{L,\mu}$. Total boundedness of $\{x \in \calX \mid L(x) \leq s, \vert \mu(x) \vert \leq t\}$ now follows by noting that it is contained in the totally bounded set $s E_{L,\mu} + \{c 1 \mid \vert c \vert \leq t \}$. 
	\end{proof}

Many interesting non-commutative examples of (M)CQMS arise from spectral triples \cite{Connes1994}, and these are the primary objects of study in the present paper. In fact, we will be focusing (almost) exclusively on the spectral triples that arise from length functions on discrete groups, for which the non-commutative metric geometry is already very complex \cite{OzawaRieffel, ChristRieffel, Klisse2026}.

	\subsection{Spectral triples and length functions on groups}\label{sec:spectral-triples}

To set the notation, we recall that a unital spectral triple $(\mathcal A,H,D)$ consists of a unital $*$-algebra faithfully represented on the Hilbert space $H$ together with a self-adjoint, possibly unbounded operator $D$ on $H$ with compact resolvent, such that $\mathcal A$ preserves $\operatorname{Dom} D$ and each commutator $[D,a]$ is bounded.
Given a spectral triple $(\mathcal A,H,D)$, we denote the $C^*$-closure of $\calA$ by $A$.  Note that one obtains a semi-norm $L_D$ on $\calA$ by setting 
\[
L_D(a)=\|[D,a]\|,
\]
which in many interesting cases, but not always \cite{Klisse2026, KyedNest}, gives rise to a CQMS. 	A spectral triple  $(\mathcal A,H,D)$ is called  a \emph{spectral metric space} if its commutator seminorm  is a Lip-norm, i.e., if $(\calA, L_D)$ is a CQMS. Amplifying to matrices, we obtain a seminorm $L_{D,n}\colon M_n(\calA) \to [0,\infty)$ by setting
\[
L_{D,n} (x):=\|[1\otimes D, x]\|
\]
and  $(\calA, H, D)$ is called a \emph{matrix spectral metric space} if $(\calA, (L_{D,n})_{n\geq 1})$  is an MCQMS; see  \cite{Kaad24}. \\
We will need to amplify unbounded operators as well. If $D\colon \operatorname{Dom}(D)\to H$ is an unbounded operator on a Hilbert space $H$ and $K$ is some auxiliary Hilbert space then we obtain an unbounded operator $D\odot 1\colon  \operatorname{Dom}(D) \odot K \to H\otimes K $ defined by $(D\odot 1)(\xi\otimes \eta)=D\xi\otimes \eta$. If $D$ is furthermore selfadjoint, then $D\odot 1$ is  essentially selfadjoint and we denote its selfadjoint closure by $D\otimes 1$.
For two spectral triples $(\calA, H_\calA, D_\calA), (\calB, H_\calB, D_\calB)$, put
	\begin{align*}
		L^{(1)}(z) = \Vert [D_\calA \odot 1, z ] \Vert, \qquad L^{(2)}(z) = \Vert [1 \odot D_\calB, z] \Vert ,
	\end{align*}
	for $z \in \calA \odot \calB$.  Note that these are indeed well defined, since for $z=\sum_{k=1}^n a_k\otimes b_k\in \calA \odot \calB $ preserves $\operatorname{Dom}(D) \odot H_{\calB} $, and on this subspace we have
\[
[D_{\calA} \odot 1, z]=\sum_{k=1}^n [D_{\calA}, a_k] \otimes b_k,
\]	
which is bounded (and similarly for $ L^{(2)}$). \\

The classical commutative example, upon which the definition of a spectral triple is modelled, arises  by choosing $\calA$ to be the algebra of smooth functions on a compact, connected spin manifold and $D$ to be the associated Dirac operator.  An interesting source of non-commutative examples is obtained by considering  length functions on discrete groups, as was already realised by Connes in his seminal paper \cite{Con:CFH}.

\begin{definition}	
	By a \emph{length function} function on a discrete group $G$ we will mean a function $\ell:G\to[0,\infty)$ with 
	\begin{enumerate}
		\item $\ell(g)=0$ if and only if $g=e$,
		\item $\ell(g^{-1})=\ell(g)$ for all $g\in G$,
		\item $\ell(gh)\leq\ell(g)+\ell(h)$ for all $g,h\in G$
	\end{enumerate}
	We say $\ell$ is \emph{proper} if the balls $B_\ell(t) := \{ g \in G\mid \ell(g) \leq t \}$ are finite for all $t \geq 0$.
\end{definition}	
    Assume that $\ell$ is proper. 
	We represent $\bbC[G]\subseteq C_r^*(G)$ on $\ell^2(G)$ through the left regular representation $\lambda_g\delta_h=\delta_{gh}$ and denote by $(\delta_g)_{g\in G}$ the standard orthonormal basis in $\ell^2(G)$. The relation $\mathcal{D}_\ell\colon \operatorname{span}\{\delta_g\mid g\in G\}\to \ell^2(G)$ defined by $\mathcal{D}_\ell\delta_g:=\ell(g)\cdot \delta_g$ is essentially selfadjoint, and its closure $D_\ell$ has compact resolvent by properness of $\ell.$
Each $\lambda_g$ preserves the domain $\operatorname{Dom} D_\ell=\left\{\xi\in \ell^2(G)\mid \sum_g\ell(g)^2|\xi(g)|^2<\infty\right\}$, and
	\begin{align}
		[D_\ell,\lambda_g]\delta_h=(\ell(gh)-\ell(h))\delta_{gh},
	\end{align}
	and hence $\|[D_\ell,\lambda_g]\|\leq\ell(g)$. Thus,  $(\bbC[G],\ell^2(G),D_\ell)$ is a unital spectral triple. Set
	\begin{align*}
		L_{\ell,n}(x):=\|[I_n\otimes D_\ell,x]\|,
		\qquad \text{and} \qquad L_\ell:=L_{\ell,1}.
	\end{align*}
We record that it is known that $(\bbC[G],\ell^2(G),D_\ell)$ is a spectral metric space if $\ell$ has bounded doubling \cite{ChristRieffel},  or  if $G$ is word-hyperbolic and $\ell$ is the word length function on $G$ \cite{OzawaRieffel} --- or more generally a length function whose induced metric is hyperbolic and weakly geodesic \cite{AustadHyp26}. A recent example by Klisse \cite{Klisse2026} shows that a word length on $(\bbZ/2 ) \wr \bbF_2$ does not yield a spectral metric space. \\
	For discrete groups $G$ and $K$ with length functions $\ell_G$ and $\ell_K$, respectively,  we obtain a natural length function on $G\times K$ by setting
	\begin{align*}
		\ell(g,k) := \ell_G(g) + \ell_K(k).
	\end{align*}
Note that the Dirac operator $D_\ell$ agrees with $D_{\ell_G} \otimes 1 + 1 \otimes D_{\ell_K}$ on the dense subspace $c_c(G)\odot c_c(K) \subseteq \ell^2(G) \tens \ell^2(K)$. Controlling  
$\Vert [D_{\ell_G} \otimes 1, z]\Vert$ and $\Vert [1 \otimes D_{\ell_K} , z]\Vert$ in terms of $\Vert [D_\ell , z]\Vert$ 
will turn out to be a crucial step in our strategy for exhibiting a CQMS structure on the product group $G \times K$ using $\ell = \ell_G + \ell_K$. We will employ a geodesic cancellation argument, see Section \ref{sec:cancelation}. \\
Lastly, note that if $\ell_G$ and $\ell_K$  are word length functions arising from  symmetric	generating sets $S_G$ and $S_K$ for $G$ and $K$, respectively, then the length function on $G\times K$ defined above is also a word length function corresponding to the generating set $\{(s,e), (e,t) \mid s\in S_G, t\in S_K\}.$

	\section{Groups and matrix compact quantum metric spaces}
	\label{sec:groups-and-MCQMS}
In this section we investigate the question of whether any CQMS is in fact an MCQMS. This question already appears in \cite{Kaad24}, in which Kaad conjectures that this is not generally the case. We first show that in the context of length functions on amenable groups, one does indeed automatically obtain an MCQMS from a CQMS,  and then proceed to show that the free group on two generators provides a counter example.

	\begin{theorem}\label{thm:amenable-CQMS-is-MCQMS}
		Let $G$ be an amenable, countable discrete group and $\ell$ a proper length function on $G$. If $(\bbC[G], L_\ell)$ is a CQMS, then $(\bbC[G], (L_{\ell,n})_{n\geq 1})$ is an MCQMS.
	\end{theorem}
	
	The proof is a modification of the argument  for \cite[Theorem 5.4]{Kaad24}.
	Write $A=C_r^*(G)$, $H=\ell^2(G)$, and $D=D_\ell$, and let $\Delta\colon\bbC[G]\to\bbC[G]\odot\bbC[G]$ be the coproduct determined by $\Delta(\lambda_g)=\lambda_g\otimes\lambda_g$.
	We shall use the isometric $*$-isomorphism
	\begin{align*}
		\iota_p\colon M_p(B(H))\otimes_{\min}B(H) \longrightarrow M_p\big(B(H)\otimes_{\min}B(H)\big),\qquad \iota_p((a_{ij})_{ij}\otimes b)=(a_{ij}\otimes b)_{ij}.
	\end{align*}
	That $\iota_p$ is indeed an  isometry follows by applying the isomorphism $M_p(E)\cong M_p(\bbC)\otimes_{\min}E$ and associativity of the minimal tensor product. We use the same notation for its restrictions to the relevant subalgebras.
	We also use the isometric map
	\begin{align*}
		\kappa_p\colon M_p(B(H))\longrightarrow M_p\big(B(H)\otimes_{\min}B(H)\big),\qquad \kappa_p((a_{ij})_{ij})=(1\otimes a_{ij})_{ij}.
	\end{align*}
	Indeed, $\kappa_p$ is the matrix amplification of the injective unital $*$-homomorphism $T\mapsto1\otimes T$ and hence isometric.
	For $\omega\in A^*$ and $z\in M_p(\bbC[G]\odot\bbC[G])$, these identifications give the identity
	\begin{equation}\label{eq:matrix-slice-identity}
		(\id_A\otimes\omega)_p(z) =(\id_{M_p(A)}\otimes\omega)\iota_p^{-1}(z).
	\end{equation}
	We note the following lemma. 
	\begin{lemma}\label{lemma:slice-lip-bound}
		For every $p\geq1$, every $\mu\in M_p(A)^*$, and every $x\in M_p(\bbC[G])$, we have
		\begin{align*}
			L_\ell\bigl((\mu\otimes\id)\iota_p^{-1}\Delta_p(x)\bigr) \leq\|\mu\|\cdot L_{\ell,p}(x).
		\end{align*}
	\end{lemma}
	\begin{proof}
		Consider the Kac--Takesaki unitary $W\in B(H\otimes H)$ given by
		\begin{align*}
			W(\delta_s\otimes\delta_t)=\delta_{t^{-1}s}\otimes\delta_t.
		\end{align*}
		It satisfies $W^*(1\otimes\lambda_g)W=\lambda_g\otimes\lambda_g$. Consequently, $\Delta$ extends to an injective unital $*$-homomorphism $A\to A\otimes_{\min}A$, given by $a\mapsto W^*(1\otimes a)W$.
		Both $W$ and $W^*$ preserve the dense subspace $\mathcal{C}:=c_c(G)\odot c_c(G)\subseteq \ell^2(G) \otimes \ell^2(G)$ and commute thereon with $1\otimes D$. For $a\in\bbC[G]$, it follows that
		\begin{equation}\label{eq:coprod-commutator}
			[1\otimes D,\Delta(a)] =W^*(1\otimes[D,a])W.
		\end{equation}
		The identity holds on $\mathcal{C}$, and hence for the bounded extensions.
		Denote by $e_{ij}$ the standard matrix units in $M_p(\bbC)$ and write $\lambda_ge_{ij}$ for the matrix in $M_p(\bbC[G])$ whose $(i,j)$-entry is $\lambda_g$ and whose other entries are zero. For $\mu\in M_p(A)^*$, we get
		\begin{align*}
			(\mu\otimes\id)\iota_p^{-1}\Delta_p(\lambda_ge_{ij}) =\mu(\lambda_ge_{ij})\lambda_g.
		\end{align*}
		We therefore have
		\begin{align*}
			[D,(\mu\otimes\id)\iota_p^{-1}\Delta_p(\lambda_ge_{ij})] &=\mu(\lambda_ge_{ij})[D,\lambda_g]\\
			&=(\mu\otimes\id) \bigl((\lambda_ge_{ij})\otimes[D,\lambda_g]\bigr)\\
			&=(\mu\otimes\id)\iota_p^{-1} \bigl([1\otimes D,\Delta(\lambda_g)]e_{ij}\bigr)
		\end{align*}    
Since these matrix units with group entries span $M_p(\bbC[G])$ we therefore obtain	
\begin{align}\label{eq:moving-mu-out}
[D,(\mu\otimes\id)\iota_p^{-1}\Delta_p(x)]=(\mu\otimes \id)\iota_p^{-1}\big(([1\otimes D, \Delta(x_{ij})])_{ij}\big)
\end{align}	
for every $x\in M_p(\bbC[G])$. 		Since any bounded functional is completely bounded with cb-norm equal to its norm \cite[Proposition~3.8]{PaulsenCBOA02}, we have $\|\mu\otimes\id\|\leq\|\mu\|$ on $M_p(A)\otimes_{\min}B(H)$. Using this bound together with the fact  that $\iota_p^{-1}$ and $\kappa_p$ are isometries and $W$ is unitary,  we now get 
		\begin{align*}
			L_\ell\bigl((\mu\otimes\id)\iota_p^{-1}\Delta_p(x)\bigr) &= \|[D, (\mu\otimes\id)\iota_p^{-1}\Delta_p(x) ] \|\\
			&=\|(\mu\otimes \id)\iota_p^{-1}\big(([1\otimes D, \Delta(x_{ij})])_{ij}\big)\| \tag{by \eqref{eq:moving-mu-out}}\\
			&\leq \|\mu\| \cdot \|([1\otimes D, \Delta(x_{ij})])_{ij}\|\\
			& =\|\mu\| \cdot \|(W^*(1\otimes [D,x_{ij}])W)_{ij}\| \tag{by \eqref{eq:coprod-commutator}} \\
			& =\|\mu\| \cdot \|((1\otimes [D,x_{ij}]))_{ij}\|\\
			& =\|\mu\| \cdot \|\kappa_p([1\otimes D, x])\|\\
			& =\|\mu\| \cdot  L_{\ell,p}(x) \qedhere
		\end{align*}

	\end{proof}

	\begin{proof}[Proof of Theorem \ref{thm:amenable-CQMS-is-MCQMS}]
		Choose a F{\o}lner sequence $(F_n)_{n\geq0}$  for $G$, with $F_0=\{e\}$. Define unit vectors and their associated vector states by
		\begin{align*}
			\xi_n=|F_n|^{-1/2}\sum_{g\in F_n}\delta_g\qquad \text{and} \qquad \chi_n(a)=\langle a\xi_n,\xi_n\rangle\quad \text{for } a\in A,
		\end{align*}
		The F{\o}lner property gives $\chi_n(\lambda_g)=|F_n\cap gF_n|/|F_n|\to1$ for every $g\in G$.
		Thus the trivial representation $\epsilon(\sum_g a_g\lambda_g):=\sum_g a_g$ is the pointwise limit
		of $\chi_n$ on $\bbC[G]$  and, by amenability,  $\epsilon$ extends  to a state on $A$.  By density, we therefore have $\lim_{n}\chi_n=\epsilon$  in the weak-$*$ topology.
		The CQMS hypothesis implies that $\delta_n:=\mk_{L_\ell}(\epsilon,\chi_n)\to0$ and by definition of $\mk_{L_\ell}$ we have
		\begin{align}\label{eq:counit-convergence}
			|(\epsilon-\chi_n)(b)|\leq\delta_n\cdot L_\ell(b)
			\qquad \text{for } b\in\bbC[G].  
		\end{align}
				Consider now the maps $\beta_n=(\id\otimes\chi_n)\Delta\colon\bbC[G]\to\bbC[G]$. They are UCP, and they have finite rank because $\beta_n(\lambda_g)=\chi_n(\lambda_g)\lambda_g$  and $\chi_n(\lambda_g)=0\quad\text{if }g\notin F_nF_n^{-1}.$
		We now show that these maps give matricial approximations. Fix $p\geq1$, $x\in M_p(\bbC[G])$, and $\mu\in M_p(A)^*$ with $\|\mu\|\leq1$. The identity $(\id\otimes\epsilon)\Delta=\id$ gives
		\begin{align*}
			|\mu(x-(\beta_n)_p(x))| &=\bigl|\mu\bigl((\id_A\otimes(\epsilon-\chi_n))_p\Delta_p(x)\bigr)\bigr|\\
			&=\bigl|\mu\bigl((\id_{M_p(A)}\otimes(\epsilon-\chi_n)) \iota_p^{-1}\Delta_p(x)\bigr)\bigr|\\
			&=\bigl|(\epsilon-\chi_n) \bigl((\mu\otimes\id)\iota_p^{-1}\Delta_p(x)\bigr)\bigr|\\
			&\leq\delta_n\cdot L_\ell \bigl((\mu\otimes\id)\iota_p^{-1}\Delta_p(x)\bigr)\\
			&\leq\delta_n \cdot L_{\ell,p}(x).
		\end{align*}
		
		The second equality is \eqref{eq:matrix-slice-identity}, the third interchanges slices on different tensor factors and the inequalities are \eqref{eq:counit-convergence} and Lemma~\ref{lemma:slice-lip-bound}, respectively.
		Since
		\begin{align*}
			\|x-(\beta_n)_p(x)\| =\sup_{\substack{\mu\in M_p(A)^*\\\|\mu\|\leq1}} |\mu(x-(\beta_n)_p(x))|,
		\end{align*}
		we obtain
		\begin{align}\label{eq:berezin-uniform-approx}
			\|x-(\beta_n)_p(x)\|\leq\delta_nL_{\ell,p}(x), 
		\end{align}
		for all $n\geq 0$ and all $p\geq 1$.
		For a given $\varepsilon>0$, choose $n$ with $\delta_n\leq\varepsilon$. Then $(\id,\beta_n)$ is a matricial $\varepsilon$-approximation in $\bbC[G]$, with the same maps for every matrix size.\\
		
		Finally, $F_0=\{e\}$ gives $\chi_0=\tau$ and $\beta_0(a)=\tau(a)1$ where $\tau$ denotes the canonical trace on $A$. 
		Applying \eqref{eq:berezin-uniform-approx} with $n=0$ yields
		\begin{align*}
			\|[x] \|\leq\|x-(\beta_0)_p(x)\| \leq \mk_{L_\ell}(\epsilon,\tau)\cdot L_{\ell,p}(x),
		\end{align*}
		and $R:= \mk_{L_\ell}(\epsilon,\tau)$ is finite because $\mk_{L_\ell}$ metrizes the weak-$*$ topology on the compact state space. This proves uniformly finite matrix diameter and completes the proof.
	\end{proof}

	Note in particular that the examples considered in \cite{ChristRieffel} have bounded doubling and are therefore amenable. As noted therein, any finitely generated group with polynomial growth has bounded doubling with respect to any word length.  We have the following immediate corollary.
	\begin{corollary}
		If $G$ is a finitely generated discrete group and $\ell$ is a proper length function with bounded doubling, then $(\bbC[G], (L_{\ell,n})_{n\geq 1})$ is an MCQMS. 
	\end{corollary}
	
	In addition to groups of polynomial growth, it is also known that word-length functions on hyperbolic groups yield compact quantum metric spaces \cite{OzawaRieffel}. Most hyperbolic groups are non-amenable, and as we shall now see the situation is quite different in this case.  Indeed, we proceed to show that $\bbF_2$ with its standard word length function does not yield an MCQMS, despite it being hyperbolic and thus yielding a  CQMS by the results in  \cite{OzawaRieffel}.
Before proceeding to the proof of this result, we	we recall the rank-one case of the projection-commutator identity in \cite[Example 6.2]{Rieffel14}.
	\begin{lemma}\label{lemma:variance}
		Let $B = B^* \in M_N(\bbC)$, and  suppose  that $v \in \bbC^N$ is a unit vector. Then the projection $p = vv^*$ satisfies
		\begin{align*}
			\Vert [B,p]\Vert^2 = \Vert Bv \Vert^2 - \vert v^* B v \vert^2.
		\end{align*}
		In particular, if $v = N^{-1/2} (1, \ldots ,1)^t$ and $B = \operatorname{diag}(b_1, \ldots, b_N)$ with $b_j \in \bbR$, then
		\begin{align*}
			\Vert [B, p ]\Vert^2 = \frac{1}{N} \sum_{j=1}^N b_j^2 - \bigg( \frac{1}{N} \sum_{j=1}^N b_j \bigg)^2. 
		\end{align*}
	\end{lemma}
	\begin{proof}
		By \cite[Example 6.2]{Rieffel14}  we have $\Vert [P,A]\Vert = \max \{ \Vert PA(1-P)\Vert, \Vert (1-P)AP \Vert \}$ for a projection $P$ and a general element $A$ in a $C^*$-algebra. For $P = p = vv^*$ and $A = B = B^*$ the two terms coincide, and we get
		\begin{align*}
			\Vert [B, p]\Vert^2 &= \sup_{\|w\|=1}\|(1-p)Bpw\|^2=\sup_{\|w\|=1}\|(1-p)B\langle w,v\rangle v \|^2=
			 \Vert (1-p)Bv \Vert^2 \\
			 &= \|Bv\|^2 - \|pBv\|^2 = \|Bv\|^2- \|\langle v,Bv  \rangle v\|^2 = \|Bv\|^2- |v^*Bv|^2	\qedhere
			 	\end{align*}
	\end{proof}
\begin{remark}\label{rem:var-interpretation}
Considering $X=\{1,\dots, N\}$ as a probability space with respect to the uniform probability measure $\mu$,  then if $b\colon X\to \bbR$ is some random variable, we have that
\[
\operatorname{Var}(b):=\mathbb{E}(b^2)- \mathbb{E}(b)^2=  \frac{1}{N} \sum_{j=1}^N b_j^2 - \bigg( \frac{1}{N} \sum_{j=1}^N b_j \bigg)^2. 
\]
The righthand side of the last formula in Lemma \ref{lemma:variance} may therefore be interpreted as the variance of a random variable, a point of view that will be important in the proof of the following result. 
\end{remark}	
	
	\begin{theorem}\label{thm:free-group-not-matrix}
		Let $G = \bbF_2$ be the free group on two generators, which we denote $a$ and $b$, and let $\ell$ be the word length for the symmetric generating set $\{a,a^{-1}, b, b^{-1}\}$. Then $(\bbC[\bbF_2], L_\ell)$ is a CQMS, but $(\bbC[\bbF_2], (L_{\ell,n})_{n\geq 1})$ is not an MCQMS. More precisely, for every operator system $Y$, every linear complete isometry $\iota \colon \bbC[\bbF_2] \to Y$, and every  finite-rank linear  map $\Phi \colon \bbC[\bbF_2] \to Y$ bounded with respect to the reduced $C^*$-norm, we have
		
		\begin{align}\label{eq:F2-failure-to-approx}
			\sup_{n\geq 1} \sup_{\substack{X\in M_n(\bbC[\bbF_2])\\ L_{\ell,n}(X)\le1}} \Vert \iota_n(X) - \Phi_n(X) \Vert \geq \frac{1}{\sqrt{3}}. 
		\end{align}
	\end{theorem}
	
	\begin{proof}
		The claim that $(\bbC[\bbF_2], L_\ell)$ is a CQMS follows from \cite[Proposition 4.3]{OzawaRieffel}. 
		Fix now an integer $R \geq 1$, and set
		\begin{align*}
			S_R = \{g\in G \mid \ell(g) = R \}.
		\end{align*}
		Since there are four choices for the first letter and three for each subsequent letter we have $N_R := \vert S_R \vert = 4 \cdot 3^{R-1}$. Then define
		\begin{align*}
			P_R = \frac{1}{N_R} [\lambda_{gh^{-1}}]_{g,h \in S_R}\in M_{N_R}(\bbC[\bbF_2]).
		\end{align*}
		Now note that the map
		\begin{align*}
			V_R \colon \ell^2(G) \to \bbC^{N_R} \otimes \ell^2(G), \quad V_R \xi = N_R^{-1/2} (\lambda_g \xi)_{g \in S_R}
		\end{align*}
		is an isometry with  $P_R = V_RV_R^*$; hence $P_R$ is a nonzero projection. We show that $L_{\ell, N_R}(P_R)\leq \sqrt{3}$. To this end, define
		\begin{align*}
			W_R := \operatorname{diag}_{g \in S_R}(\lambda_g) \in M_{N_R}(\bbC[\bbF_2]), \quad v_R := N_R^{-1/2} (1, \ldots , 1)^t\in \bbC^{N_R} \quad \text{and} \quad p_R := v_R v_R^*\in M_{N_R}(\bbC). 
		\end{align*}
		Then $W_R^* P_R W_R = \big[ \frac{1}{N_R} 1 \big]_{g,h \in S_R}  = p_R \otimes 1$. 
		Note that $W_R$ and $W_R^*$ preserve the domain of $I_{N_R} \otimes D$. On the summand $\bbC^{N_R} \otimes \bbC \delta_x$ the operator $W_R^* (I_{N_R} \otimes D) W_R$ acts as $B_x := \operatorname{diag}_{g \in S_R} \ell(gx)$. Using that $W_R$ is unitary, Lemma \ref{lemma:variance} then gives 
		\begin{align}
			L_{\ell, N_R}(P_R)^2 &= \Vert [I_{N_R}\otimes D_\ell, P_R]\Vert^2 \notag \\
			&=\Vert [W_R^*(I_{N_R}\otimes D_\ell)W_R, p_R \otimes 1]\Vert^2 \notag \\
			&= \sup_{x \in \bbF_2} \Vert [B_x, p_R] \Vert^2\\
			& = \sup_{x \in \bbF_2} \operatorname{Var} (\ell( - x)) \label{eq:LNR-variance-bound},
		\end{align}
		where we in the last equality consider the function $S_R\ni g\mapsto \ell(gx)\in \bbR$ as a random variable on $S_R$ endowed with the uniform probability measure; see e.g.~Remark \ref{rem:var-interpretation}. \\
		
		For $x\in \bbF_2$ and $g\in S_R$,  denote by $C(g,x)$ the number of cancellation pairs in the reduced product $gx$. Then
		\begin{align*}
			\ell(gx) = R + \ell(x) - 2C(g,x).
		\end{align*}
		We now fix an $x\in \bbF_2$ and consider the random variable $C\colon S_R\to \bbR$ given by $C(g):=C(g,x)$.
		If the reduced word for $x$ is $t_1\cdots t_m$, the event $C \geq k$, for $1 \leq k \leq \min(R,m)$, prescribes the suffix $t_k^{-1}\cdots t_1^{-1}$ to $g$. There are $3^{R-k}$ such words, so
		\begin{align*}
			\mathbb{P}(C \geq k) = \frac{3^{R-k}}{4\cdot 3^{R-1}}= \frac{1}{4 \cdot 3^{k-1}}. 
		\end{align*}
		For larger $k$ the probability is zero. We now note that $C^2=\sum_{k=1}^{\infty} (2k-1)1_{C\geq k}$; indeed, if $C(g)=j$ then it holds that  $j^2 = \sum_{k=1}^j (2k-1)$. We thus obtain
		\begin{align*}
			\mathbb{E}[C^2] = \mathbb{E}\Big[\sum_{k=1}^\infty (2k-1) 1_{C \geq k}\Big] = \sum_{k=1}^\infty (2k-1) \mathbb P(C \geq k) \leq \frac{1}{4} \sum_{k=1}^\infty (2k-1) 3^{-(k-1)} = \frac{3}{4}. 
		\end{align*}
		Now, from $\ell(gx) = R + \ell(x) -2C(g)$ we obtain $\ell(-x) - \mathbb E(\ell(-x)) = -2(C - \mathbb E(C))$. Thus 
		\begin{align*}
			\operatorname{Var}(\ell(-x)) &= \mathbb{E}\big[\big(\ell(-x)- \mathbb{E}[\ell(-x)]\big)^2 \big]= 4 \mathbb E[(C - \mathbb E[C])^2] \\
			&= 4 \operatorname{Var}(C) =  4\big(\mathbb E[C^2] - (\mathbb E[C])^2 \big)\leq 4 \mathbb E[C^2] \leq 3.
		\end{align*}
		The assertion $L_{\ell, N_R}(P_R) \leq \sqrt{3}$ now follows from \eqref{eq:LNR-variance-bound}.\\
		
		Consider now a bounded finite rank map $\Phi\colon \bbC[\bbF_2] \to Y$ for some operator system $Y$.  We now aim to show that 
		\[
		\lim_{R\to \infty} \|\Phi_{N_R}(P_R)\|=0.
		\]
		 An application of \cite[Lemma 1.4]{HaagerupNCMAP}  shows that
		\begin{align}\label{eq:Haagerup-ineq}
			\big \Vert \sum_{g \in S_R} c_g \lambda_g \big\Vert \leq (R+1) \bigg( \sum_{g \in S_R} \vert c_g \vert^2 \bigg)^{1/2}.
		\end{align}
		for any family $(c_g)_{g\in S_R}$ of complex numbers. 
		For  a unit vector $\alpha\in \bbC^{N_R}$, set
		\begin{align*}
			U_\alpha := \sum_{g \in S_R} \overline{\alpha_g} \lambda_g
		\end{align*}
		Then for any $f\in A^*$ and any two unit vectors $\alpha, \beta\in \bbC^{N_R}$, we have
		\begin{align*}
		\langle f_{N_R}(P_R)\beta, \alpha \rangle =	\alpha^* f_{N_R}(P_R) \beta = \frac{1}{N_R} \sum_{g,h} \overline{\alpha_g} \beta_h f(\lambda_{gh^{-1}} ) = \frac{1}{N_R} f(U_\alpha U_\beta^*).
		\end{align*}
		Applying \eqref{eq:Haagerup-ineq} to $U_\alpha$ and $U_\beta^*$ we obtain
		\begin{align}\label{eq:free-RD-upper-bound}
		\|f_{N_R}(P_R)\|&=\sup_{\|\alpha\|=\|\beta\|=1} |\langle f_{N_R}(P_R)\beta, \alpha \rangle| =	\sup_{\|\alpha\|=\|\beta\|=1} \frac{1}{N_R}\cdot |f(U_\alpha U_\beta^*)|\notag\\
		&\leq \frac{1}{N_R} \sup_{\|\alpha\|=\|\beta\|=1} \|f\|\cdot\|U_\alpha\|\cdot \|U_\beta\|
				 \leq \Vert f \Vert \cdot \frac{(R+1)^2}{N_R}. 
		\end{align}
		Now choose a basis $y_1, \ldots , y_s$ of $\operatorname{Ran} \Phi$. We may then write $\Phi(a) = \sum_{j=1}^sf_j (a) y_j$ for some $f_j \in A^*$. Then \eqref{eq:free-RD-upper-bound} yields
		\begin{align*}
			\Vert \Phi_{N_R}(P_R) \Vert \leq \sum_{j =1}^s \Vert (f_j)_{N_R} (P_R) \Vert \Vert y_j \Vert \leq \bigg( \sum_{j=1}^s \Vert f_j \Vert \Vert y_j \Vert \bigg) \frac{(R+1)^2}{4 \cdot 3^{R-1}} \longrightarrow 0. 
		\end{align*}
		Finally, set $X_R = P_R/\sqrt{3}$. Then $L_{\ell, N_R} (X_R) \leq 1$, and $\iota$ being a complete isometry implies
		\begin{align*}
			\Vert \iota_{N_R} (X_R) - \Phi_{N_R}(X_R) \Vert \geq \frac{1}{\sqrt{3}} (1- \Vert \Phi_{N_R}(P_R) \Vert ). 
		\end{align*}
		Letting $R \to \infty$ proves \eqref{eq:F2-failure-to-approx}. In particular, the finite rank UCP maps needed to establish an MCQMS cannot give an approximation with error constant smaller than $1/\sqrt{3}$.   
	\end{proof}

	\section{Cartesian products}
	\label{sec:Cartesian-products}
	Our next aim is to show how one may use the theory of matrix compact quantum metric spaces to expand on the class of groups with length functions that yield compact quantum metric spaces. As mentioned already,  the only structural results in this direction cover the classes of hyperbolic groups \cite{OzawaRieffel} and of groups of polynomial growth \cite{ChristRieffel}. We first develop some general machinery that can be used to verify that a  certain length function on a cartesian product of two groups yields a CQMS, and then proceed to verify that the relevant assumptions  are met in a number of instances. \\

	\subsection{A tensor product criterion}
	\label{tensor:onefactor}

	We first recall some standard properties of the minimal tensor product.
	Let $A$, $B$, and $C$ be unital $C^*$-algebras. For every $w\in C\otimes_{\min}B$, we have
	\begin{align}\label{eq:norm-achieving-UCP}
		\Vert w\Vert =\sup\big\{  \Vert(\mathrm{id}_C\otimes\psi)(w)\Vert \mid  n\geq 1, \psi\in\mathrm{UCP}(B,M_n(\mathbb C)) \big\}
	\end{align}
which may be seen by representing $B$ faithfully on a Hilbert space and considering the UCP maps arising from compressions to finite-dimensional subspaces.

	By injectivity of the minimal tensor product, a unital complete isometry $j\colon A\to C$ induces a complete isometry $j\otimes\mathrm{id}_B:A\otimes_{\min}B\to C\otimes_{\min}B$,  see \cite[Theorem~4.6]{KPTT11}.
	Finally, every bounded linear functional $f\in C^*$ satisfies $\Vert f\Vert_{\mathrm{cb}}=\Vert f\Vert $ by \cite[Proposition~3.8]{PaulsenCBOA02}.
	The completely bounded mapping property of the minimal tensor product (see e.g.~\cite{Pisier2003}) therefore gives
	\begin{align*}
		\Vert (f\otimes\mathrm{id}_B)(w)\Vert \leq\Vert f\Vert \,\Vert w\Vert,
		\qquad w\in C\otimes_{\min}B.
	\end{align*}

	\begin{proposition}
		\label{prop:one-leg-approx}
		Let $(\mathcal X,L_\calX)$ and $(\mathcal Y,L_\calY)$ be CQMSs 
		and let $M$ be a Lipschitz seminorm on
		$\mathcal X\odot\mathcal Y$. Suppose that there exist
		$c_{\mathcal X},c_{\mathcal Y}>0$, such that
		\begin{align*}
			L_{\mathcal X}((\id\otimes g)z)&\leq c_{\mathcal X}\Vert g\Vert M(z) &&(g\in \mathcal{\mathcal Y } ^*),\\
			L_{\mathcal Y}((f\otimes\id)z)&\leq c_{\mathcal Y}\Vert f\Vert M(z) &&(f\in \mathcal{X}^*).
		\end{align*}
		Suppose, moreover, that for every $\varepsilon>0$ there exist a unital
		$C^*$-algebra $C_\varepsilon$, a unital complete isometry
		$\iota_\varepsilon \colon \mathcal X\to C_\varepsilon$, and a bounded linear map
		$\Phi_\varepsilon\colon \mathcal X\to C_\varepsilon$ of finite rank satisfying
		\begin{equation}\label{eq:one-leg-epsilon-approx}
			\Vert ((\iota_\varepsilon-\Phi_\varepsilon)\otimes\id)z\Vert 
			\leq\varepsilon M(z),\qquad z\in\mathcal  X\odot\mathcal Y.
		\end{equation}
		Then $(\mathcal X\odot\mathcal Y,M)$ is a CQMS.
	\end{proposition}
Here we consider $\mathcal X\odot\mathcal Y$ with the operator system structure induced by the embedding into the minimal tensor product of their enveloping $C^*$-algebras.

	\begin{proof}
		Choose states $\mu \in S(\mathcal X)$ and $\nu \in S(\mathcal Y)$. By Proposition \ref{prop:totally-bounded-criterion} it suffices to show that
		\begin{align*}
			E_{M, \mu\otimes \nu} = \{ z \in \calX \odot \calY \mid M(z) \leq 1, \  (\mu \otimes \nu)(z) = 0 \} 
		\end{align*}
		is totally bounded in norm. Fix $\varepsilon>0$ and the maps $\iota_\varepsilon, \Phi_\varepsilon \colon \mathcal{X} \to C_\varepsilon$ ensured by the assumptions. 
		Choose a basis $e_1, \ldots , e_r$ of $\operatorname{Ran}\Phi_\varepsilon$. We may then write $\Phi_\varepsilon(x) = \sum_{i=1}^r f_i(x) e_i$ with $f_i \in \mathcal{X}^*$. Defining $b_i(z) = (f_i \otimes \id )z$, we then have
		\begin{align*}
			(\Phi_\varepsilon \otimes \id)z = \sum_i e_i \otimes b_i(z), \qquad \text{and} \qquad L_\calY(b_i(z)) \leq c_\calY \Vert f_i \Vert, 
		\end{align*}
		for $z \in E_{M, \mu\otimes \nu}$. 
		For $z \in E_{M, \mu\otimes \nu} $, we now define $a_z:=(\id\otimes\nu)z$. Then $L_{\mathcal{X}}(a_z)\leq c_\calX$ and $\mu(a_z)=0$. 
		Setting $E_{L_{\mathcal{X}}, \mu} := \{ x \in \calX \mid L_{\mathcal{X}}(x) \leq 1, \mu(x) = 0 \}$ and $R_{\mathcal{X}} := \sup_{x \in E_{L_{\mathcal{X}}, \mu}}\Vert x \Vert$, we have that $a_z/c_{\mathcal{X}}\in E_{L_{\mathcal{X}}, \mu}$ and 
		\begin{align*}
			\vert \nu(b_i(z)) \vert = \vert f_i(a_z) \vert \leq c_\calX R_\calX \Vert f_i \Vert. 
		\end{align*}
		Since 
		\[
		\{b_i(z) \mid z \in E_{M, \mu \otimes \nu} \} \subseteq c_\calY \Vert f_i \Vert E_{L_\calY, \nu} + \{\lambda 1_\calY \mid \vert \lambda \vert \leq c_\calX R_\calX \Vert f_i \Vert \},
		\]
		 we deduce that the former set is totally bounded for $i=1, \ldots ,r$.  Hence, $(\Phi_\varepsilon\otimes \id)(E_{M, \mu \otimes \nu})$ is totally bounded as well. 
				
		We may therefore choose $z_1 ,\ldots z_N \in E_{M, \mu\otimes \nu}$ such that $(\Phi_\varepsilon\otimes \id) z_1, \dots, (\Phi_\varepsilon\otimes \id)z_N$form an $\varepsilon$-net in $(\Phi_\varepsilon\otimes \id)(E_{M, \mu \otimes \nu})$. Since $\iota_\varepsilon\otimes \id$ is isometric, for each $z \in E_{M, \mu \otimes \nu}$  and a suitable $i\in \{1,\dots, N\}$ we have
		\begin{align*}
			\Vert z - z_i \Vert &=\Vert (\iota_\varepsilon \otimes \id) z - (\iota_\varepsilon \otimes \id) z_i\Vert \\
			&\leq \Vert (\iota_\varepsilon \otimes \id) z - (\Phi_\varepsilon\otimes \id)z \Vert + \Vert (\Phi_\varepsilon\otimes \id)z - (\Phi_\varepsilon\otimes \id) z_i \Vert + \Vert (\Phi_\varepsilon\otimes \id)z_i - (\iota_\varepsilon \otimes \id)z_i \Vert \\
			&< 3 \varepsilon
		\end{align*}
		by using \eqref{eq:one-leg-epsilon-approx}. It follows that $E_{M, \mu \otimes \nu}$ is totally bounded. 
	\end{proof}
	
We now specialise to spectral triples, which are the main case of interest in the present paper.	For the statement, recall that a spectral triple $(\calA, H_A, D_A)$ is called a (matrix) spectral metric space if the commutator seminorm yields a (M)CQMS; see Section \ref{sec:spectral-triples} for details
	
	\begin{theorem}\label{thm:one-MCQMS-suffices}
		Suppose $(\calA, H_A, D_A)$ is a matrix spectral metric space and $(\calB, H_B, D_B)$ is a spectral metric space. For $z = \sum_k a_k \otimes b_k \in \calA\odot \calB$, set
		\begin{align*}
			P(z) = \sum_k [D_A, a_k] \otimes b_k \qquad \text{and} \qquad Q(z) = \sum_k a_k \otimes [D_B, b_k].
		\end{align*}
		If $M$ is a Lipschitz seminorm on $\calA \odot \calB$ for which there exist $c_A, c_B > 0$ such that
		\begin{align}\label{eq:one-MCQMS-suffices-int-eq}
			\Vert P(z) \Vert \leq c_A M(z) \qquad \text{and} \qquad \Vert Q(z) \Vert \leq c_B M(z),
		\end{align}
		for all $z \in \calA \odot \calB$, then $(\calA \odot \calB , M)$ is a CQMS. The analogous assertion holds with the two factors interchanged. 
	\end{theorem}
	
	\begin{proof}
		We verify the hypotheses of Proposition \ref{prop:one-leg-approx}. Denote by $A$ and $B$ the respective $C^*$-completions, and put $L_{D_A,n} (x) = \Vert [I_n \otimes D_A , x ] \Vert$. Let $\delta > 0$ be given, and use that $(\calA, H_A, D_A)$ is a matrix spectral metric space to obtain a $C^*$-algebra $C$ along with maps $\iota, \Phi \colon \calA \to C$ such that
		\begin{align}\label{eq:iota-phi-eq-one-MCQMS-suffices}
			\Vert \iota_n (x) - \Phi_n (x) \Vert \leq \delta \cdot L_{D_A,n} (x)
		\end{align}
		for all $n \geq 1$ and all $x \in M_n(\calA)$. Here $\iota$ is a complete isometry and $\Phi$ is UCP with finite-dimensional range; in particular both maps extend to $A$.\\
		
		Fix now an $n\in \mathbb{N}$, a UCP map $\psi \colon B \to M_n(\bbC)$ as well as a $z\in \calA \odot \calB$, and set $x_\psi = (\id \otimes \psi)z$. Under the canonical matrix identification $\calA \odot M_n(\bbC) \cong M_n(\calA)\cong M_n(\bbC) \odot \calA$ we then have
		\begin{align}\label{eq:moving-psi-out}
			[1_n \otimes D_A , x_\psi] = (\id_{B(H_A)} \otimes \psi)P(z). 
		\end{align}
		Setting $W = ((\iota - \Phi) \otimes \id_B)(z) \in C \otimes_{\rm min} B$, we now obtain 
		\begin{align*}
			\Vert (\id_C \otimes \psi)(W) \Vert = \Vert \iota_n(x_\psi) - \Phi_n(x_\psi) \Vert \leq \delta \cdot L_{D_A,n} (x_\psi) \leq \delta \cdot \Vert P(z) \Vert,
		\end{align*}
		where the first inequality follows from \eqref{eq:iota-phi-eq-one-MCQMS-suffices} and the last inequality from \eqref{eq:moving-psi-out}. 
		Since this holds for all $n\in \mathbb{N}$ and all $\psi \in \operatorname{UCP}(B, M_n(\bbC))$, \eqref{eq:norm-achieving-UCP} and \eqref{eq:one-MCQMS-suffices-int-eq} imply
		\begin{align*}
			\Vert ((\iota-\Phi)\otimes \id)z \Vert \leq \delta \cdot \Vert P(z) \Vert \leq \delta c_A \cdot M(z). 
		\end{align*}
		Given $\varepsilon>0$, choosing $\delta = \varepsilon/c_A$ we see that \eqref{eq:one-leg-epsilon-approx} is satisfied. Moreover, 
		\begin{align*}
			L_{D_A}((\id \otimes g) z) \leq \Vert g \Vert \Vert P(z) \Vert \qquad \text{and} \qquad L_{D_B}((f \otimes \id)z) \leq \Vert f \Vert \Vert Q(z) \Vert,
		\end{align*}
		and so Proposition \ref{prop:one-leg-approx} applies and shows that  $(\calA \odot \calB , M)$ is a CQMS. It is clear that the argument is symmetric in the choice of matrix factor. 
	\end{proof}

	\subsection{Cancellation in one factor}\label{sec:cancelation}
	We now return to the setting of groups and length functions to provide a criterion ensuring that the assumptions in Theorem \ref{thm:one-MCQMS-suffices} are satisfied. 
	Suppose therefore  $G$ and $K$ are countable discrete groups equipped with proper length functions $\ell_G$ and $\ell_K$, and consider the length function $\ell_{G \times K} := \ell_G + \ell_K$ on $G\times K$. Suppressing the natural isomorphism   $\ell^2(G \times K) \cong \ell^2(G) \otimes \ell^2(K)$, we then have
	\[
	D_{\ell_{G\times K}}= D_{\ell_G} \otimes 1 + 1 \otimes D_{\ell_K}, 
	\]
as unbounded operators on, say, $c_c(G) \odot c_c(K)$. For notational convenience, we set
	
	\begin{align*}
		D_1 := D_{\ell_G} \otimes 1, \qquad D_2 := 1 \otimes D_{\ell_K}, \quad \text{ and } \quad  D:=D_{\ell_{G\times K}}.
	\end{align*}
	For $x \in M_p(\bbC[G \times K])$, we furthermore set
	\begin{align*}
		T^1(x) := [1_p \otimes D_1, x], \qquad T^2(x) := [1_p \otimes D_2, x] \quad \text{and} \quad  T (x) := T^1 (x) + T^2 (x),
	\end{align*}
	which are all extendable to bounded operators. We make no notational difference between these densely defined operators and their bounded extensions.
		The strategy in the sequel will be to average certain sequences of unitary conjugates of $T(x)$ which fix $T^1(x)$ and make $T^2(x)$ converge strongly to zero. This is formalised by means of the following definition. 
	
	\begin{definition}
		A \emph{finite cancellation datum} for $(K, \ell_K)$ consists of a finite nonempty set $I$, numbers $w_i \geq 0$ with $\sum_i w_i = 1$, and sequences $(g_i(n))_n \subseteq K$ indexed by $i\in I$, such that, for every fixed $m,s \in K$,
		\begin{align}\label{eq:geodesic-datum-cond}
			\sum_{i\in I}w_i\bigl( \ell_K(msg_i(n))-\ell_K(sg_i(n)) \bigr)
			\underset{n\to \infty}{\longrightarrow} 0.
		\end{align}
		It is \emph{eventually exact} if this expression is zero for all sufficiently large $n$, with the threshold allowed to depend on $m$ and $s$.
	\end{definition}
	
	\begin{lemma}\label{lemma:cancellation-datum-implication}
		If $K$  admits a finite cancelation datum
		then for every $p \geq 1$ and every $x \in M_p (\bbC[G \times K])$, it holds that
		\begin{align*}
			\Vert T^1(x) \Vert \leq \Vert T(x) \Vert \qquad \text{and} \qquad \Vert T^2(x) \Vert \leq 2 \Vert T  (x) \Vert. 
		\end{align*}
	\end{lemma}
	
	For the proof, we will need the right translations as well. For $t \in K$, we therefore let $\rho_t \delta_s = \delta_{st}$ be the right translation, which is a unitary operator. It commutes with every left translation $\lambda_m$, and just like for left translation, every $\rho_t$ preserves the domain of $D_{\ell_K}$. 
	\begin{proof}
		Write $x = \sum_{m \in F} x_m \otimes \lambda_m$ with $x_m \in M_p (\bbC[G])$ for some finite set $F \subseteq K$. To simplify notation, we will  denote $1_p\otimes 1_{\ell^2(G)}\tens \rho_t$ by $V_t$ in the following. 
		 Since $V_t$ operator commutes with $x$ and $D_1$, it will fix $T^1(x)$ when conjugating, that is $V_t^* T^1(x) V_t = T^1(x)$. For $\eta \in \bbC^p \otimes \ell^2(G)$ and $s \in K$, we have
		\begin{align*}
			V_t^* T^2(x) V_t (\eta \otimes \delta_s) = \sum_{m \in F} (\ell_K(mst) - \ell_K(st)) x_m \eta \otimes \delta_{ms}. 
		\end{align*}
		Notably, the rightmost tensor-leg coordinate, $ms$, does not depend on $t$. For each $n\in \mathbb{N}$, we now define a contraction on $\mathbb{B}(\bbC^p\tens \ell^2(G)\tens \ell^2(K))$ by
		\begin{align*}
			\Psi_n (S) = \sum_i w_i V_{g_i(n)}^* S V_{g_i(n)}.
		\end{align*}
		For $\xi = \sum_{s \in E}\xi_s \otimes \delta_s\in \bbC^p\tens \ell^2(G)\tens \ell^2(K)$ with $E$ finite, we have
		\begin{align*}
			\Vert \Psi_n (T^2(x)) \xi \Vert \leq \sum_{s\in E, m\in F} \Vert x_m \Vert \cdot \Vert \xi_s \Vert \cdot \big| \sum_i w_i (\ell_K(msg_i(n)) - \ell_K (sg_i(n))) \big| \longrightarrow 0. 
		\end{align*}
		Such vectors are dense, and $\Vert \Psi_n (T^2(x)) \Vert \leq \Vert T^2 (x) \Vert$, so we have strong convergence everywhere. Consequently, $\Psi_n (T(x)) \to_n T^1(x)$ strongly. Contractivity then gives $\|T^1(x) \xi \Vert \leq \Vert T(x) \Vert \Vert \xi \Vert$ for every $\xi$, which proves the first inequality in the lemma statement. The second then follows by writing $T^2(x) = T(x) - T^1(x)$. 
	\end{proof}
\begin{corollary}\label{cor:consequence-for-grps}
If $(\bbC[G], \ell^2(G), D_{\ell_G})$ is a spectral metric space and if $(\bbC[K], \ell^2(K), D_{\ell_K})$ is a matrix spectral metric space and if $K$ admits a finite cancelation datum, then $(\bbC[G\times K], \ell^2(G\times K), D_{\ell_{G\times K}})$ is a spectral metric space. 
\end{corollary}	
\begin{proof}
By Lemma \ref{lemma:cancellation-datum-implication}, the assumptions in Theorem \ref{thm:one-MCQMS-suffices} are satisfied and the result follows. 
\end{proof}

	To illustrate cancellation, we include the following easy lemma.
	\begin{lemma}\label{lemma:finite-group-cancellation}
		
		For any  length function on a finite group $K$, one obtains an exact cancellation datum by setting $I=K$, $w_k=\tfrac{1}{|K|}$ and $g_k(n)=k$ for all $n\in \mathbb{N}$. 
	\end{lemma}
	\begin{proof}
		For fixed $m,s$ the maps $t \mapsto mst$ and $t \mapsto st$ both permute $K$, and so the sums of their lengths are equal. 
	\end{proof}
	For abelian and infinite dihedral groups, we construct geodesic rays along which the length defects eventually stabilize. Averaging over the two directions and a finite set of residues then gives cancellation.
	The construction is related to Rieffel's horofunction methods \cite{Rieffel02}.
	For an equally weighted family it suffices to prove, for each fixed $h$,
	\begin{align}\label{eq:def-An}
		A_n(h):=\sum_i\bigl(\ell_K(hg_i(n))-\ell_K(g_i(n))\bigr)=0 \quad\text{eventually}.
	\end{align}
	Indeed, the sum in \eqref{eq:geodesic-datum-cond} is $\bigl(A_n(ms)-A_n(s)\bigr)/|I|$.\\
	
	We now proceed to show that certain natural groups admit a finite cancelation datum.

	\subsection{Finitely generated abelian groups}
	We first show that finitely generated abelian groups have a finite cancelation datum. 
	\begin{proposition}\label{prop:fga-cancellation}
		Every finite symmetric word length on a finitely generated abelian group $K$ has an eventually exact cancellation datum. If $K$ is infinite, the sequences can be chosen as samples of geodesic rays. 
	\end{proposition}
	\begin{proof}
		Since Lemma \ref{lemma:finite-group-cancellation} covers the case for finite $K$, we assume $K$ is infinite in the sequel. 
		We will write $K$ additively and denote its rank by $d$ and its torsion subgroup by $\operatorname{Tor}(K)$. Suppose $S$ is the given generating set, and consider the quotient map
		\begin{align*}
			\pi \colon K \to K/ \operatorname{Tor}(K) \cong \bbZ^d \subseteq \bbR^d 
		\end{align*}
		for some $d \geq 1$. The set $P :=\operatorname{conv}(\pi(S))$ is then a polytope of dimension $d$,
		 which is centrally symmetric and it contains $0$ in its interior. Recall that a facet of $P$ is a face of dimension $d-1$. Choose a facet $F$ together with a supporting linear functional $\sigma \colon \bbR^d \to \bbR$; that is $\sigma(p) \leq 1$  for $x\in P$ and $F = P \cap \{\xi \in \bbR^d \mid \sigma(\xi) = 1 \}$; see e.g.~\cite{Ziegler1995} for background on polytopes.\\
		Since $0$ is in the interior of $P$, the maximum value for $\sigma$ is positive. 
				 To ease notation we will also write $\sigma(w) = \sigma(\pi(w))$ for $w \in K$. Now define
		\begin{align*}
			V := \{v\in S \mid \pi(v) \in F \}, \quad q := |V|, \quad  z:= \sum_{v \in V} v \quad \text{ and } \quad \Lambda:= \langle V \rangle\leq K. 
		\end{align*}
		Note that if $w = s_1 + \cdots + s_k$ with $s_i \in S$, then $\sigma(w) = \sum_i \sigma(s_i) \leq k$. Taking the minimum over all representing words, we obtain $\ell_K(w) \geq \sigma(w)$. Since  $z$ is defined as the sum of $q$ letters from $V$,  $nz$ can be represented by $nq$ such letters. Since $\sigma(nz) = nq$, we deduce $\ell_K (nz) = nq$ for all $n\geq 0$. As $S = -S$, we deduce $\ell_K(-w) = \ell_K(w) $ for every $w \in K$. If $n< 0$, this means $\ell_K(nz) = \ell_K((-n)z) = (-n) q$, and so
		\begin{align}\label{eq:fga-length-expression}
			\ell_K (nz) = |n|q \qquad \text{for all } n\in \bbZ. 
		\end{align}
		Let $\gamma \colon \bbN_0 \to K$ be the path obtained by successively adding the letters of this repeated listing of elements of $V$, starting at $0$; i.e.~writing $V=\{s_1,\dots, s_q\}$ the path is defined as : 
	\[	
	\gamma(0)=0, \quad \gamma(1)=s_1, \quad \gamma(2)=s_1+s_2,  \cdots\cdots,  \gamma(q)=z, \quad \gamma(q+1)=z+s_1, \cdots
	\]	
		If $r < t$, the subword from $\gamma(r)$ to $\gamma(t)$ has $t-r$ letters, each of which takes the value $1$ under $\sigma$. Hence
		\begin{align*}
			t-r = \sigma(\gamma(t) - \gamma(r)) \leq \ell_K(\gamma(t) - \gamma(r)) \leq t-r.
		\end{align*}
		It follows that $\gamma$ is geodesic, and $\gamma(nq) = nz$.\\
		
		We proceed to show that $K/\Lambda$ is finite. Note that $F = \operatorname{conv}(\pi(V))$ and that it has dimension $d-1$, and hence $\pi(V)$ contains $d$ affinely independent points $u_1, \ldots u_d$. 
		By affine independence, a relation $\sum_{j=1}^dt_ju_j=0$ with real coefficients satisfying $\sum_{j=1}^dt_j=0$ forces every $t_j$ to vanish. In fact, these points are also linearly independent. We may see this by applying $\sigma$ to any relation $\sum_{j=1}^dt_ju_j=0$ and noting that since $\sigma(u_j) = 1$, this gives $\sum_{j=1}^dt_j=0$. Affine independence then forces $t_1 = \cdots = t_d= 0$.
		Consequently, $\pi(\Lambda)$ contains $d$ linearly independent integer vectors, $w_1,\dots, w_d$. The matrix $A$ whose columns consist of $w_1,\dots, w_d$ is therefore invertible and its inverse is given by $A^{-1}=(\det(A))^{-1}\operatorname{Adj}(A)$, where the adjugate matrix $\operatorname{Adj}(A)$ has integer entries. We therefore have that
		\[
			\det(A)\cdot \bbZ^d =A\cdot \operatorname{Adj}(A) \cdot \bbZ^d\subseteq A\bbZ^d\subseteq \pi(\Lambda),
		\]
		and it follows that $\pi(\Lambda)$  has finite index in $\bbZ^d$. 
		The surjection $K/\Lambda\to\mathbb Z^d/\pi(\Lambda)$ has kernel $(\Lambda+\operatorname{Tor}(K))/\Lambda$ which, by the second isomorphism theorem, is isomorphic to $\operatorname{Tor}(K)/\Lambda\cap \operatorname{Tor}(K)$, and thus finite. Finiteness of $K/\Lambda$ now follows from noting that both kernel and range are finite. \\
		
		For fixed $w \in K$, consider the integer-valued sequence
		\begin{align*}
			b_n (w) = \ell_K(nz + w) - nq.
		\end{align*}
		Then $\ell_K ((n+1) z + w) \leq \ell_K (nz+w) + q$ implies $b_{n+1}(w) \leq b_n(w)$, so the sequence is non-increasing. 
		Using \eqref{eq:fga-length-expression} we also have $nq = \ell_K(nz) \leq \ell_K(nz + w) + \ell_K(w)$, and so we deduce 
		$-\ell_K(w)\leq b_n(w)\leq\ell_K(w)$, forcing the sequence to eventually stabilize at a value, which we denote by $b(w)$. 
				We will need to understand how $b$ changes under translations with elements from $\Lambda$. For $v \in V$, we have $\ell_K(v) = 1$. Moreover, as $z- v = \sum_{u\in V\setminus \{v\}} u$ is a sum of $q-1$ generators, we have $\ell_K(z-v) \leq q-1$. Writing $(n+1)z + w = (nz + w + v) + (z-v)$, we use the triangle inequality to calculate
		\begin{align*}
			\ell_K (nz + w + v) &\leq \ell_K(nz + w) +1, \\
			\ell_K((n+1)z + w) &\leq \ell_K(nz + w + v) + q - 1.
		\end{align*}
		Thus
		\begin{align*}
			b_{n+1}(w) +1 \leq b_n(w+v) \leq b_n (w) +1.
		\end{align*}
		As $b_n$ eventually stabilizes to $b$, we see that for sufficiently large $n$, both the outer terms equal $b(w) +1$, and so $b(w+v) = b(w) +1$. We also note that this yields $b(w) - 1 = b(w-v +v) -1 = b(w-v)$. Every $a \in \Lambda$ can be written as $a = \sum_{v\in V} k_v v$ for $k_v \in \bbZ$. By iterating these identities we get
		\begin{align*}
			b(w+a) = b(w) + \sum_{v\in V} k_v = b(w) + \sum_{v\in V} k_v \sigma(v) = b(w) + \sigma\Big(\sum_{v\in V} k_v v\Big) = b(w) + \sigma(a). 
		\end{align*}
		It follows that $b(w+a) - \sigma(w+a) = b(w) - \sigma(w)$, and so $f([w]) := b(w) - \sigma(w)$ is a well-defined function on $K/\Lambda$. Thus, for all sufficiently large $n$ (the threshold depending on $w$) we have
		\begin{align}\label{eq:fga-length-2}
			\ell_K(nz + w ) = nq + \sigma(w) + f([w]).
		\end{align}

		With the analysis above at our disposal, we now obtain the  claimed cancellation datum as follows: Choose as set of representatives, $R$, for $K/\Lambda$, put $M = |R|$, and define $2M$ sequences as follows:
		\begin{align*}
			g^+_r (n) &:= \gamma(nq) +r =  nz +r,    &\text{for } r\in R,\\
			 g^-_r(n)  &:= -\gamma(nq) -r = -nz -r, &\text{for } r\in R.
		\end{align*}
		Fix $h \in K$, and choose $n_0\in \bbN$ such that   \eqref{eq:fga-length-2}  holds for all $n\geq n_0$ and all $w$ in the finite set $\{r, r \pm h\mid r\in R\}$. Then 
		\begin{align*}
			\sum_{r\in R}\bigl(\ell_K(h+nz+r)-\ell_K(nz+r)\bigr) &=\sum_{r\in R}\bigl(\sigma(h)+f([h+r])-f([r])\bigr) =M\sigma(h),\\
			\sum_{r\in R}\bigl(\ell_K(h-nz-r)-\ell_K(-nz-r)\bigr) &=\sum_{r\in R}\bigl(-\sigma(h)+f([r-h])-f([r])\bigr) =-M\sigma(h).
		\end{align*}
		For the second line, symmetry gives $\ell_K(h-nz-r)=\ell_K(nz+r-h)$ and $\ell_K(-nz-r)=\ell_K(nz+r)$.
		In both lines the $f$-terms sum to zero because $[r]\mapsto[r\pm h]$ permutes the finite quotient $K/\Lambda$. The two sums therefore cancel, giving $A_n(h)=0$, where $A_n$ is given as in \eqref{eq:def-An}. 
		The weighted sum in \eqref{eq:geodesic-datum-cond} is then
		\begin{align*}
			\frac{A_n(m+s)-A_n(s)}{2M}=0.
		\end{align*}
		Thus, the weights $1/(2M)$ give the required eventually exact datum. 
	\end{proof}
	
	\subsection{The infinite dihedral group}
	Next we show that also the infinite dihedral group $D_\infty$ admits a finite cancelation datum. 
	\begin{proposition}\label{prop:geodesic-cancellation-inf-dihedral}
	
		Let $K$ be the infinite dihedral group
		\begin{align*}
			K = \langle a,b \mid b^2 = e, bab = a^{-1} \rangle,
		\end{align*}
		and let $\ell_K$ be the word length for any finite symmetric generating set $S$. Then there is an eventually exact cancellation datum whose sequences are samples of geodesic rays. 
	\end{proposition}
	
	\begin{proof}
		Write elements of $K$ as $a^j b^\varepsilon$ for $j \in \bbZ$ and $\varepsilon \in \{0,1\}$, and note that
		\begin{align*}
			(a^jb^\delta)(a^kb^\epsilon)  =a^{j+(-1)^\delta k}b^{\delta+\epsilon},
		\end{align*}
		where the exponent of $b$ is understood to be modulo $2$. Now define
		\begin{align*}
			T&=\{t\in \bbZ :a^t\in S\}  \quad \ \ \text{ and } \quad   P=\max(\{|t|:t\in T\}\cup\{0\}),\\
			 Q&=\{r\in \bbZ: a^rb\in S\} \quad \text{ and } \quad \Delta=\max Q-\min Q.
		\end{align*}
		Note that $T = -T$. Also, $Q$ is nonempty since $K$ is non-abelian, and $\nu := \max \{P ,\Delta/2 \} >0$ because otherwise $S$ would generate a group of order at most $2$. Now choose positive integers $q, L$ by
		\begin{align*}
			(q,L) = \begin{cases}
				(P,1), & \text{if } P\geq \Delta/2 \\
				(\Delta,2), &\text{if } P < \Delta/2.
			\end{cases}  
		\end{align*}
		Then $q= \nu L$. 
		Note that if $(q,L) = (P,1)$, then $a^q = a^P \in S$ as $T = -T$. If $(q,L) = (\Delta, 2)$, then $(a^{\max Q}b)(a^{\min Q}b) = a^{\Delta} = a^q$ and $(a^{\max Q}b)(a^{\min Q}b)$ has length  2. Thus, in either case, there is a word of length $L$ representing $a^q$. 
				We first prove 
		\begin{align*}
			\ell_K(a^{nq}) = |n| L, \qquad \text{for all } n\in \bbZ. 
		\end{align*}
		To this end, put $c := (\max Q + \min Q) /2$ and define $X\colon K\to \bbR$  by $X(a^k b^\varepsilon) := k - c \varepsilon$.  We now wish to analyse how much the value of $X$ changes upon multiplying the input by a generator.  We first note that
		\begin{align*}
			X(a^k b^\varepsilon a^t) - X(a^k b^\varepsilon) &= (-1)^\varepsilon t, \\
			X(a^k b^\varepsilon a^rb) - X(a^k b^\varepsilon) &= (-1)^\varepsilon(r-c). 
		\end{align*}
		If $a^t\in S$ then $|t|\leq P\leq \nu$  and if $a^rb\in S$ then $|r-c| \leq \Delta/2\leq \nu$.
		So,  upon multiplying a given element  $a^k b^\varepsilon\in K$ by a generator from the right, the absolute change in the value of $X$ is at most $\nu$.
		 Also, $X(e) = 0$, so telescoping along a word of length $k$ representing $h$ gives $\vert X(h) \vert \leq k \nu$. Thus $\ell_K(h) \geq \vert X(h) \vert /\nu$.  If $h = a^{nq}$ for some $n\in \bbZ$, this reads $\ell_K (a^{nq}) \geq |n|q/\nu = |n|L$, which proves the lower bound.  But, as we saw above $a^q$ can be represented by a word of length $L$ so $a^{nq}$ can be expressed as a product of $|n|L$ generators,  which provides the upper bound. \\

		For $\sigma \in \{\pm1 \}$, $0\leq r <  q$ and $\varepsilon \in \{0,1\}$, set
		\begin{align*}
			g_{\sigma,r,\varepsilon}(n) := a^{\sigma n q + r} b^\varepsilon.
		\end{align*}
		We note that
		\begin{align*}
			g_{\sigma,r,\varepsilon}(n)^{-1} g_{\sigma,r,\varepsilon}(k) = a^{(-1)^\varepsilon \sigma(k-n)q}.
		\end{align*}
		The length of this element is $|n-k|L$ by what we noted above. In other words, the distance between $g_{\sigma,r,\varepsilon}(n)$ and $g_{\sigma,r,\varepsilon}(k)$ in the induced left-invariant word metric is $|n-k|L$.
		Repeating a shortest word for $a^{(-1)^\epsilon\sigma q}$ from the initial point $a^rb^\epsilon$ interpolates these samples to a geodesic ray. Indeed, the interpolating path has length $|k-n|L$ between any two samples. 
		A shorter path between two intermediate vertices, together with the remaining portions of the interpolating path, would contradict the distance formula for two enclosing samples.
		For fixed $(\sigma,r,\epsilon)$, the integers
		\begin{align*}
			f_n=\ell_K(a^{\sigma nq+r}b^\epsilon)-nL
		\end{align*}
		are non-increasing, since left multiplication by $a^{\sigma q}$ increases the length by at most $L$. Using subadditivity and the reverse triangle inequality we obtain
		\begin{align*}
			nL - \ell_K(a^r b^\varepsilon) \leq \ell_K (a^{\sigma n q + r} b^\varepsilon) \leq nL + \ell_K(a^r b^\varepsilon). 
		\end{align*}
		It follows that $f_n$ is bounded and integer-valued, and must therefore stabilize. We may thus  find an $N_0$ and constants $C_{\sigma,r\varepsilon} \in \bbZ$ such that
		\begin{align}\label{eq:stabilized-geodesics}
			\ell_K(g_{\sigma,r\varepsilon}(n) ) = nL + C_{\sigma,r,\varepsilon} \qquad \text{for all } n\geq N_0,  
		\end{align}
		for all $4q$ indices $(\sigma, r,\varepsilon)$. \\

		With this, we may now finally verify cancellation. Fix $h = a^j b^\delta$. For every  $r\in\{0, \dots, q-1\}$, the there exist unique $t_r\in \bbZ$ and $r'\in\{0,\dots, q-1\}$ such that
		\begin{align*}
			j+(-1)^\delta r=qt_r+r'.
		\end{align*}
		We now define $\sigma':=(-1)^\delta\sigma$ and $\epsilon':=\delta+\epsilon \, \,(\operatorname{mod} 2)$. Then
		\begin{align*}
			h g_{\sigma, r , \varepsilon}(n) = g_{\sigma' , r', \varepsilon'} (n + \sigma' t_r) .
		\end{align*}
		Note that $t_r$ is independent of $\sigma$ and $\varepsilon$. For $n \geq N_0 + \max_r \vert t_r \vert$, both $n$ and $n + \sigma' t_r$ are at least $N_0$, and so \eqref{eq:stabilized-geodesics} yields
		\begin{align*}
			\ell_K(h g_{\sigma, r, \varepsilon}(n)) - \ell_K (g_{\sigma, r, \varepsilon}(n)) = \sigma' t_r L + C_{\sigma' ,r', \varepsilon'} - C_{\sigma, r, \varepsilon}. 
		\end{align*}
		The map of indices $(\sigma, r, \epsilon)$ to $(\sigma',r',\epsilon')$ is a permutation: from its image one recovers $\sigma=(-1)^\delta\sigma'$, $\epsilon=\epsilon'-\delta$ modulo two, and the unique $r\in\{0,\ldots,q-1\}$ congruent to $(-1)^\delta(r'-j)$ modulo $q$.
		Consequently the sums of the two families of constants agree.
		For fixed $r,\epsilon$, summing over $\sigma\in\{1,-1\}$ cancels the carry terms because $t_r$ is fixed and $\sigma'$ takes both signs.
		Thus
		\begin{align}\label{eq:dihedral-sum}
			\sum_{\sigma,r,\epsilon} \bigl(\ell_K(hg_{\sigma,r,\epsilon}(n)) -\ell_K(g_{\sigma,r,\epsilon}(n))\bigr) = 0
		\end{align}
		for every sufficiently large $n$.  This is $A_n(h)=0$, so the weights
		$1/(4q)$ give \eqref{eq:geodesic-datum-cond}.

	\end{proof}
	
	\subsection{Products with hyperbolic groups}
	In this section we show how the analysis carried out above leads to new examples of length functions on groups which give rise to compact quantum metric spaces. We first summarize the upshot of the previous section as follows: 
	
	\begin{proposition}
		\label{prop:bounds-established}
		Let $G$ be any group with  a proper length function $\ell_G$.  Let $K$ be either
		\begin{itemize}
		\item[-] a finite group with any length function, or
		\item[-] a finitely generated abelian group with any finite symmetric word length, or 
		\item[-] the infinite dihedral group with any finite symmetric word length. 
		 \end{itemize}
  Then the bounds of Lemma \ref{lemma:cancellation-datum-implication} hold at every matrix level, with constants $1$ and $2$ independent of the group, generating set, and matrix size.
	\end{proposition}
	
	\begin{proof}
		This is immediate by Lemma \ref{lemma:cancellation-datum-implication}, Lemma \ref{lemma:finite-group-cancellation}, Proposition \ref{prop:fga-cancellation} and Proposition \ref{prop:geodesic-cancellation-inf-dihedral}.
	\end{proof}
	
	\begin{theorem}
		\label{thm:cartesian-products}
		Let $G$ be a finitely generated group equipped with a word-length function $\ell_G$ for which $(\mathbb C[G],L_{\ell_G})$ is a CQMS, and
		let  $K_1,\dots, K_m$  be a family or groups, each of which is either  finite,  finitely generated abelian or isomorphic to $D_\infty$. Equipping each $K_j$ with a word length function $\ell_{K_j}$ and
		\begin{align*}
			P:=G\times K_1\times\cdots\times K_m \text{ with }\ell(g,k_1,\ldots,k_m)=\ell_G(g)+\sum_{j=1}^m\ell_{K_j}(k_j),
		\end{align*}
		the pair $(\mathbb C[P],L_\ell)$ is a CQMS. 
	\end{theorem}
	
	\begin{proof}
		Each $K_j$ is amenable and yields a CQMS by \cite{ChristRieffel}. Theorem \ref{thm:amenable-CQMS-is-MCQMS} therefore makes every $K_j$ an MCQMS. For $m=1$, the result therefore follows by Corollary \ref{cor:consequence-for-grps}. The result now follows by induction on $m$.
	\end{proof}
	
	\begin{corollary}\label{cor:cartesian-prod-hyperbolic}
		If $H$ is word-hyperbolic and $K_1,\ldots,K_m$ are as in Theorem \ref{thm:cartesian-products}, then
		$H\times K_1\times\cdots\times K_m$, with length function coming from the union of the chosen generating sets for the factors, yields a CQMS.
	\end{corollary}
	
	\begin{remark}
		Corollary \ref{cor:cartesian-prod-hyperbolic} covers groups of the form $G = H \times F \times D_\infty^k \times \bbZ^d $ for $H$ word hyperbolic and $F$ finite, and with the length function as described in the result. In particular it covers the case $d = 0$, $F$ a finite Coxeter group, and $H$ a non-elementary hyperbolic Coxeter group, and $k\geq 1$. As such, we have the first examples of word length functions on non-amenable,  non-hyperbolic Coxeter groups with CQMS structures for these word-lengths. Indeed, these groups contain $H$ and a copy of $\bbZ^2$, and are therefore neither amenable nor word-hyperbolic. They are thus in particular not covered by the results of \cite{ChristRieffel} or \cite{OzawaRieffel}.
	\end{remark}


\begin{thebibliography}{KPTT11}

\bibitem[AKK22]{AKK:Podcon}
Konrad Aguilar, Jens Kaad, and David Kyed.
\newblock The {P}odle\'{s} spheres converge to the sphere.
\newblock {\em Comm. Math. Phys.}, 392(3):1029--1061, 2022.

\bibitem[Aus26]{AustadHyp26}
Are Austad.
\newblock Compact quantum metric spaces from weakly geodesic length functions
  on hyperbolic groups, 2026.

\bibitem[Con89]{Con:CFH}
Alain Connes.
\newblock Compact metric spaces, {F}redholm modules, and hyperfiniteness.
\newblock {\em Ergodic Theory Dynam. Syst.}, 9(2):207--220, 1989.

\bibitem[Con94]{Connes1994}
Alain Connes.
\newblock {\em Noncommutative Geometry}.
\newblock Academic Press, San Diego, 1994.

\bibitem[CR17]{ChristRieffel}
Michael Christ and Marc~A. Rieffel.
\newblock Nilpotent group {${\rm C}^*$}-algebras as compact quantum metric
  spaces.
\newblock {\em Canad. Math. Bull.}, 60(1):77--94, 2017.

\bibitem[Cv21]{walter-connes:truncations}
Alain {Connes} and Walter~D. {van Suijlekom}.
\newblock {Spectral truncations in noncommutative geometry and operator
  systems}.
\newblock {\em {Commun. Math. Phys.}}, 383(3):2021--2067, 2021.

\bibitem[Haa79]{HaagerupNCMAP}
Uffe Haagerup.
\newblock An example of a nonnuclear {$C\sp{\ast} $}-algebra, which has the
  metric approximation property.
\newblock {\em Invent. Math.}, 50(3):279--293, 1978/79.

\bibitem[Kaa24]{Kaad24}
Jens Kaad.
\newblock External products of spectral metric spaces.
\newblock {\em Pure Appl. Funct. Anal.}, 9(5):1175--1211, 2024.

\bibitem[KK21]{KaadKyed2021}
Jens Kaad and David Kyed.
\newblock Dynamics of compact quantum metric spaces.
\newblock {\em Ergodic Theory and Dynamical Systems}, 41(7):2069--2109, 2021.

\bibitem[Kli26]{Klisse2026}
Mario Klisse.
\newblock Word-length spectral triples of
  {$(\mathbb{Z}/2\mathbb{Z})\wr\mathbb{F}_{d}$} are not metric, 2026.

\bibitem[KN24]{KyedNest}
David Kyed and Ryszard Nest.
\newblock Finiteness of metrics on state spaces.
\newblock {\em Bull. Lond. Math. Soc.}, 56(1):288--295, 2024.

\bibitem[KPTT11]{KPTT11}
Ali Kavruk, Vern~I. Paulsen, Ivan~G. Todorov, and Mark Tomforde.
\newblock Tensor products of operator systems.
\newblock {\em J. Funct. Anal.}, 261(2):267--299, 2011.

\bibitem[Lei25]{Leimbach2025}
Malte Leimbach.
\newblock Convergence of {Peter--Weyl} truncations of compact quantum groups.
\newblock {\em Journal of Noncommutative Geometry}, 2025.
\newblock Published online first.

\bibitem[OR05]{OzawaRieffel}
Narutaka Ozawa and Marc~A. Rieffel.
\newblock Hyperbolic group {$C^*$}-algebras and free-product {$C^*$}-algebras
  as compact quantum metric spaces.
\newblock {\em Canad. J. Math.}, 57(5):1056--1079, 2005.


\bibitem[Pau02b]{PaulsenCBOA02}
Vern Paulsen.
\newblock {\em Completely bounded maps and operator algebras}, volume~78 of
  {\em Cambridge Studies in Advanced Mathematics}.
\newblock Cambridge University Press, Cambridge, 2002.

\bibitem[Pis03]{Pisier2003}
Gilles Pisier.
\newblock {\em Introduction to Operator Space Theory}, volume 294 of {\em
  London Mathematical Society Lecture Note Series}.
\newblock Cambridge University Press, 2003.

\bibitem[Rie02]{Rieffel02}
Marc~A Rieffel.
\newblock Group {C*}-algebras as compact quantum metric spaces.
\newblock {\em Documenta Mathematica}, 7:605--651, 2002.

\bibitem[Rie14]{Rieffel14}
Marc~A. Rieffel.
\newblock Standard deviation is a strongly {L}eibniz seminorm.
\newblock {\em New York J. Math.}, 20:35--56, 2014.

\bibitem[{van}21]{walter:GH-convergence}
Walter~D. {van Suijlekom}.
\newblock {Gromov-Hausdorff convergence of state spaces for spectral
  truncations}.
\newblock {\em {J. Geom. Phys.}}, 162:11, 2021.
\newblock Id/No 104075.

\bibitem[Zie95]{Ziegler1995}
G{\"u}nter~M. Ziegler.
\newblock {\em Lectures on Polytopes}, volume 152 of {\em Graduate Texts in
  Mathematics}.
\newblock Springer, New York, 1995.

\end{thebibliography}

\end{document}